\documentclass[a4paper,12pt,twoside]{article}
\usepackage[latin1]{inputenc}
\usepackage[T1]{fontenc}
\usepackage{amssymb}
\usepackage{dochier}
\usepackage[final]{lpretex}
\usepackage{title}
\usepackage{a4}
\usepackage{amscd}
\usepackage{caption}
\input xy
\xyoption{all}

\EndOfDump

\def\DHpartformat#1{(#1) }
\def\DHrefpart#1{(\DHRefpart{#1})}
\LBsetstyleof{partno}{arabic}

\renewcommand{\labelenumi}{(\arabic{enumi})}

\def\i {\item}

\let\define\def

  \def\C {{\mathbb C}}
  
\def\GG {{\mathbb G}} \def\H{{\mathbb H}}  
  
  \def\P {{\mathbb P}} 
\def\Q {{\mathbb Q}} \def\R {{\mathbb R}}
 
  \def\X {{\mathbb X}}
\def\Z {{\mathbb Z}} 

\define \n {\mathbb N}
\define \z {\mathbb Z}
\define \q {\mathbb Q}
\define \PP {\mathbb P}

\def\sA {{\Cal A}}  \def\sC {{\Cal C}}
\def\sD {{\Cal D}} \def\sE {{\Cal E}} \def\sF {{\Cal F}}
  
  \def\sL {{\Cal L}}
\def\sM {{\Cal M}}  \def\sO {{\Cal O}}
 
\def\sS {{\Cal S}}  
  \def\sX {{\Cal X}}

\define \cN {\Cal N}
\define \cf {\Cal F}
\define \cg {\Cal G}
\define \cE {\Cal E}
\define \ce {\Cal E}
\define \cc {\Cal C}
\define \cV {\Cal V}
\define \cA {\Cal A}
\define \cK {\Cal K}
\define \cO {\Cal O}
\define \cF {\Cal F}
\define \cn {\Cal N}
\define \cI {\Cal I}
\define \sP {\Cal P}
\define \sEll {\Cal{Ell}}
\define \sJE {\Cal{JE}}

\define \sGJE {\Cal{GJE}}
\define \sHyp {\mathcal{HYP}}
\def\a {\alpha} \def\b {\beta} \def\g {\gamma}  
   
\def\s {\sigma}

\define \x {\xi}
\define \y {\eta}
\define \G {\Gamma}
\define \r {\rho}
\define \w {\omega}

\def \trho {\tilde {\rho}}

\define \tH {\widetilde H}
\define \tG {\widetilde{\Gamma}}
\define \tW {\widetilde W}
\define \tF {\widetilde F}
\define \tm {\tilde m}
\define \St {\widetilde S}
\define \Xt {\widetilde X}
\define \tS {\widetilde S}
\define \tpsi {\tilde \psi}
\define \tL {\widetilde L}
\define \tE {\widetilde E}
\define \tl {\tilde l}
\define \tA {\widetilde A}
\define \tom {\tilde\omega}
\define \tT {\widetilde T}
\define \tB {\widetilde B}
\define \tf {\tilde f}
\define \tsA {\widetilde{\sA}}

\define \tM {\widetilde M}
\define \tpsi {\widetilde{\psi}}
\define \trho {\widetilde{\rho}}
\define \tR {\widetilde R}
\define \tp {\widetilde p}
\define \tq {\widetilde q}
\define \tc {\tilde c}
\define \tsF {\widetilde {\sF}}
\define \tsM {\widetilde {\sM}}
\define \tii {\tilde i}
\define \tx {\tilde x}
\define \tg {\tilde g}
\define \tw {\tilde w}
\define \tz {\tilde z}
\define \ta {\widetilde\alpha}

\define \tmu{\widetilde{\mu}}

\define \bD {\overline{D}}
\define \bG {\overline{G}}
\define \bI {\overline{I}}
\define \bK {\overline{K}}

\define\sje{\sJE}
\define\bsell{\bsEll}
\define \bV {\overline{V}}
\define \bX {\overline{X}}
\define \bY {\overline{Y}}
\define \btau {\overline{\tau}}
\define\bareta{\overline{\eta}}

\def\pd {\partial}

\def \Dx1 {\frac{\pd}{{\pd} x_1}}
\def \Dy1 {\frac{\pd}{{\pd} y_1}}
\def \Dz1 {\frac{\pd}{{\pd} z_1}}
\def \Dx2 {\frac{\pd}{{\pd} x_2}}
\def \Dy2 {\frac{\pd}{{\pd} y_2}}
\def \Dz2 {\frac{\pd}{{\pd} z_2}}

\def\q {\quad} 

\def\Mapdiagr#1{\nearrow\rlap{$\lower 5pt\vbox{{\hbox{$\mkern
-15mu\scriptstyle#1$}}}$}} 

\def\Mapdiagl#1{\llap{$\lower 5pt\vbox{{\hbox{$\scriptstyle#1\mkern
-15mu$}}}$}\searrow} 

\def\Mapswr#1{\swarrow\rlap{$\lower 5pt\vbox{{\hbox{$\mkern
-15mu\scriptstyle#1$}}}$}}              

\def\Mapnwl#1{\nwarrow\rlap{$\lower 5pt\vbox{{\hbox{$\mkern
-15mu\scriptstyle#1$}}}$}}

\def\i.e#1#2#3{\mathrel{\smash{\mathop{#2}\limits^{#1}_{#3}}}}

\def \inj {\hookrightarrow}

\define \Rhook {\hookrightarrow}

\def \half {\raise1pt\hbox{$\scriptstyle
        \frac{1}{2}\displaystyle$}}

\def \x{{\sl X}\llap{$\mkern -2mu {\scriptstyle -}$}}

\def \Hom {\operatorname{Hom}}

\define \Kod {\operatorname{Kod}}
\define \dimension {\operatorname{dim}}
\define \codim {\operatorname{codim}}
\define \contr {\operatorname{contr}}
\define \rk {\operatorname{rank}}
\define \Im {\operatorname {Im}}
\define \Mor {\operatorname{Mor}}
\define \Cl {\operatorname{Cl}}
\define \Hilb {\operatorname{Hilb}}
\define \degree {\operatorname{deg}}
\define \mult {\operatorname{mult}}
\define \Aut {\operatorname{Aut}}
\define \NS {\operatorname{NS}}
\define \Gal {\operatorname{Gal}}
\define \ch {\operatorname{char}}
\define \Jac {\operatorname{Jac}}
\define \Km {\operatorname{Km}}
\define \Sec {\operatorname{Sec}}
\define \Stab {\operatorname{Stab}}
\define \Br {\operatorname{Br}}
\define \Inv {\operatorname {Inv}}
\define \tr {\operatorname{tr}}
\define \Frob {\operatorname{Frob}}
\define \Symn {\operatorname{Symm}^n}
\define \Ev {\sE^\vee}
\define \ordp {\operatorname{ord}_p}
\define \Supp {\operatorname{Supp}}
\define \Ann {\operatorname{Ann}}
\define \disc {\operatorname{disc}}
\define \lie {\operatorname{lie}}
\define \embdim {\operatorname{embdim}}

\def\Ram{\operatorname{Ram}}

\define\HS{\operatorname{HS}}

\def\bsEll{\overline{\sE\ell\ell}}

\define\nbd{neighbourhood }

\def\hod#1#2#3#4{\ensuremath{ if#30 H^{#2}({#1},{\cal O}_{#1}) \else 
 H^{#2}(#1,\Omega^{#3}\if\relax{#4}\relax_{#1}\else _{#1/#4}\fi)\fi}}

\begin{document}
\title[Torelli for elliptic surfaces]
{Generic Torelli with denominators for  
 elliptic surfaces}
\author{N. I. Shepherd-Barron}
\address{King's College,
Strand,
London WC2R 2LS,
U.K.}
\email{Nicholas.Shepherd-Barron@kcl.ac.uk}
\maketitle
\begin{abstract}
We show that a very general Jacobian
elliptic surface is determined by its polarized
$\Q$-Hodge structure, subject to various constraints
on the irregularity and the geometric genus.
\end{abstract}
AMS classification: 14C34, 32G20.
\begin{section}{Introduction}\label{intro}
  We consider the stack $\sje_{h,q}$ of Jacobian elliptic surfaces
  $f:X\to C$ over $\C$
  of geometric genus $h$ and irregularity $q$,
  provided that $8h> 10(q-1)$ and
  $h\ge q+3$.
  If all the singular fibres are of Kodaira type $I_1$
  then the datum of such a surface is equivalent to the datum of
  a classifying morphism $\Phi:C\to \bsell$ from $C$ to the stack
  $\bsell$ of generalized elliptic curves.

  In \cite {SB} we showed that then the derivative of the period map
  associated  to the primitive cohomology $H^2_{prim}(X)$,
  which is defined as the orthogonal complement
  of a fibre and the zero section,
  determines the base curve $C$, the ramification divisor
  $Z=\Ram_\Phi$ in $C$  of $\Phi$ and the
  copy of $C$ in $\P^{h-1}$ that arises as the image of $X$
  under the linear system
  $\vert K_X\vert = f^*\vert K_C+\Phi^*M\vert$.
  (Here $M$ is the bundle  of weight one modular forms on $\bsell$
  and, under our assumptions, $\vert K_X\vert$ has no
  base points, so that the infinitesimal Torelli
  theorem holds for these surfaces.)
  We went on to show
  that if also $q\ge 2$ and the pair $(C,Z)$ is generic
  (as we assume for the rest of this paper)
  then these data determine $\Phi$ modulo the automorphism group
  $\GG_m$ of $\bsell$. We then proved, without the assumption
  that $q\ge 2$, that the generic Torelli theorem holds for these
  surfaces.
  In this paper we consider such surfaces under the 
  assumptions that
  \begin{eqnarray}\label{assume}
    8h>10(q-1),\ h\ge q+3\ \textrm{and}\ q\ge 2&{}&
  \end{eqnarray}
  \noindent and we show that
  a very general such surface $f:X\to C$
  is determined by its polarized $\Q$-Hodge structure.
  To do so it is enough, given the results that we have just
  described,
  to assume that the base curve $C$ and the line bundle
  $L=\Phi^*M$ on $C$ are fixed and to
  consider only surfaces in the $\GG_m$-orbit
  of the point defined by $f:X\to C$ in the closed substack
  $\sje_{C,L}$ of the stack
  $\sje=\sje_{h,q}$ of these surfaces.

    Recall (\cite{Mir89}, p. 28) that
    $X$ is embedded as a divisor
    in the $\P^2$-bundle
    $$\pi:\P=\P_{C,L}=\P(\sO\oplus L^{2}
    \oplus L^{3})\to C$$
    and is linearly equivalent to the line bundle
    $\pi^*(L^6)\otimes\sO_{\P}(3)$.
    It is defined by an
    equation
    $$Y^2Z = 4X^3 +g_4 XZ^2 + g_6 Z^3$$
    where $g_n\in H^0(C, L^{\otimes n})$.
    If $g_6=0$ then $X$ is a \emph{Gauss surface}
    and if $g_4=0$ then $X$ is an \emph{Eisenstein surface}
    because in that case $X$ has an action by
    $\Z/4$ or $\Z/6$, respectively. In the Eisenstein
    case the action of $\langle\a\rangle=\Z/6$
    is given by $\a^*(X,Y,Z)=(\zeta_6^2 X, \zeta_6^3Y, Z)$.
      
        The closure of the $\GG_m$-orbit
        through $X$ is the pencil $\{X_{(\lambda,\mu)}\}_{(\lambda,\mu)\in\P^1}$
        given by the equation
        $$Y^2Z = 4X^3 +\lambda g_4 XZ^2 + \mu g_6 Z^3,$$
        where $\lambda,\mu$
        are homogeneous co-ordinates on the
        copy of $\P^1$ that is the closure of $\GG_m$.
        We call this the \emph{Gauss--Eisenstein pencil
          through $X$} and regard it as a canonically defined
        pencil through
        $X$ that is generated by the Eisenstein surface $X_6$
        defined by $\lambda=0$ and the Gauss
        surface $X_4$ defined by
        $\mu=0$. The surfaces $X_6$ and $X_4$ are determined
        by $X$, since they are the fixed points in the closure of
        the $\GG_m$-orbit through $X$.
        When $g_6$ has isolated zeroes then $X_6$
        is smooth and when $g_6$ has a double
        zero then $X_6$ has a singularity of type $A_2$.

        Until now we have only considered general surfaces in $\sje$.
        The condition that the surface be general excludes Gauss and Eisenstein
        surfaces.
        To include them, we consider the stack $\sje^{RDP}$
        whose geometric points are relative canonical models
        of Jacobian elliptic surfaces. (These are the surfaces
        with RDPs that arise from contraction of the vertical
        $(-2)$-curves that are disjoint from the zero section.
        For such surfaces there
        might be no classifying morphism $\Phi:C\to\bsell$,
        but there is still a zero section and its conormal bundle,
        so it makes sense to speak of
        the closed substack $\sje^{RDP}_{C,L}$ of $\sje^{RDP}$.)
        The Gauss--Eisenstein pencils are in these
        stacks $\sje^{RDP}_{C,L}$. If $\sD$ is the relevant
        period domain (in the context of this paper it classifies Hodge structures on
        $H^2_{prim}(Y,\Z)$ when $Y$ is a Jacobian elliptic surface
        of geometric genus $h$ and irregularity $q$)
        and $\G$ the relevant arithmetic group
        then the period map
        exists as a morphism
        $$per:[\sje^{RDP}]\to [\sD/\G]$$
        of geometric quotients; here, and for the
        rest of this paper, $[\sX]$ denotes
        the geometric quotient of the Deligne--Mumford stack
        $\sX$ in either of the algebraic or the analytic contexts.
\end{section}
\begin{section}{Eisenstein surfaces}
  Fix both a curve $C$ (of
  genus $q\ge 0$) and a line bundle $L$ on it of degree $h+1-q>0$.
  We consider Jacobian elliptic surfaces $X\to C$ such that
  $p_g(X)=h$ and
  $L$ is the conormal bundle of the zero section.
  Note that $L$
  also the pull-back $\Phi^*M$ of the bundle of
  weight $1$ modular forms.
  \begin{lemma} If $L$ is generic then $6L$
    is very ample.
    \begin{proof} This follows from the assumptions that
      $8h> 10(q-1)$ and $h\ge q+3$.
    \end{proof}
  \end{lemma}

  We shall consider more particularly the corresponding family
  $\{X_t\}$ of Eisenstein surfaces defined by $g_4=0$; the parameter
  space is the linear system $\vert 6L\vert$. We put
  $V=H^2_{prim}(X_{\bareta}, \Q)$,
  where $X_{\bareta}$ is the geometric generic
  member of this family, and let $G$ 
  denote the monodromy group acting on $V$.
  The action of $\g =\a^2$ on $X_{\bareta}$
  makes
  $V$ into a $\Q[T]/(T^3-1)$-module and then
  $G$ is a group of $\Q[T]/(T^3-1)$-linear transformations
  of this module.
  In this section we re-write that part (4.4) of \cite{Del80}
  that is concerned with monodromy in even
  fibre dimension in this context.

  Take a general pencil $\Pi$ in $\vert 6L\vert$.
  Since $6L$ is very ample
  this gives a family
  $$a:\sX\to\Pi\cong \P^1$$ whose total space $\sX$ is a
  smooth blow-up of $\P_{C,L}$ and
  whose geometric generic member
  $X_{\bareta}$ is smooth. Moreover, again because
  $6L$ is very ample, each of the
  finitely many singular fibres
  has a single $A_2$-singularity.
  Note that $V$ is identified with the orthogonal
  complement of the image of $H^2(\sX,\Q)$
  in $H^2(X_{\bareta},\Q)$.

  Let $F$ denote a local Milnor fibre in this family,
  so that $H^2(F,\Q)$ is the $\Q$-module
  spanned by roots $\delta_1,\delta_2$ such that
  $\delta_i^2=-2$ and $\delta_1.\delta_2 = 1$.
  Then
  the local monodromy operator $\s$ on $H^2(F,\Q)$
  is a Coxeter element in the Weyl group $W(A_2)$,
  so is of order $3$.
  We shall refer to these roots $\delta_i$
  as the \emph{basic cycles} and to the Coxeter elements
  $\s$ as \emph{basic transformations}.

  By letting $\Pi$
  move in the Grassmannian of lines in $\vert 6L\vert$
  we see, as in \cite{SGA7 II} XVIII,
  that the basic transformations
  form a single conjugacy class in
  $G$ and they
  generate $G$.

  \begin{lemma}\label{vanish} $V^G=0$.
    \begin{proof}
      Since the fibres of $\sX\to \Pi$ all have only du Val
      singularities, they satisfy the hypotheses of
      Th{\'e}or{\`e}me 1.1 of \cite{SGA7 II} XVIII, so that
      the Leray spectral sequence
      $$E_2^{pq}=H^p(\Pi, R^qa_*\Q) \Rightarrow
      H^{p+q}(\sX,\Q)$$
      degenerates at $E_2$. From this it follows
      that
      $H^2_{prim}(X_{\bareta},\Q)^G$ consists of classes that
      come from the threefold $\sX$
      and the lemma is proved.
    \end{proof}
  \end{lemma}

  \begin{lemma}\label{Milnor}
    The action of $\g$ on $H^2(F,\Q)$
    is non-trivial.
    \begin{proof} $F$ is the affine surface
      defined by the equation $xy+z^3=1$
      and $\g$ acts as $\g^*(x,y,z)=(x,y,\zeta_6^4 z)$.
      So the fixed locus of $\g$ in $F$ is $\C -\{0\}$,
      whose Euler characteristic is zero.
      We conclude by the Lefschetz fixed point theorem.
    \end{proof}
  \end{lemma}

  It follows that the action of
  $\g$ on $H^2(F,\Q)$ is given by
  $\g=\s^{\pm}$.

  Recall that $\Z/3=\langle\g\rangle$ has just two irreducible
  representations over $\Q$, namely, the trivial
  representation
  $V_0=\Q$ and the $2$-dimensional representation
  $V_1=\Q(\zeta_3)$, on which $\g$ acts as multiplication by
  $\zeta_3$.
  In particular, $H^2(F,\Q)$ is isomorphic to $V_1$; in other words,
  $H^2(F,\Q)$ is a $1$-dimensional $\Q(\zeta_3)$-vector space and as
  such is generated
  by a basic cycle.

  The next lemma is well known.
  \begin{lemma} If $A$ and $B$ are groups, $k$ is a field
    and $W$ is a representation of $A\times B$ defined
    over $k$ on which the action of $A$ is completely
    reducible, then there is a decomposition
    $$W=\oplus_j\left( V_j\otimes_k U_j\right)$$
    of $k[A\times B]$-modules,
    where $V_j$ runs over the irreducible representations
    of $A$ defined over $k$ and
    $U_j$ is the representation of $B$ given by
    $U_j=\Hom_A(V_j,W)$.
    \noproof
  \end{lemma}

  Taking $k=\Q$,
  $A=\langle\g\rangle$, $B=G$ and $W=V$ gives
  $V=(V_0\otimes U_0) \oplus (V_1\otimes U_1)$
  for some representations $U_0, U_1$ of $G$.
  \begin{lemma} $V_0\otimes U_0=V^\g=0$.
    \begin{proof} That $V_0\otimes U_0= V^\g$
      follows from the definition of $U_0$ and $V_0$, while $V^\g$
      is also the primitive cohomology
      of the geometric quotient $[X/\g]$.
      This is a $\P^1$-bundle over $C$, so that
      its primitive cohomology is zero.
    \end{proof}
  \end{lemma}

  \begin{corollary} $V$ is naturally a representation
    of $G$ defined over $\Q(\zeta_3)$.
    \begin{proof} We already know that $V$ is a $\Q[T]/(T^3-1)$-module
      and we have just shown that the eigenspace belonging
      to $T=1$ is zero.
    \end{proof}
  \end{corollary}
  The intersection pairing $b:V\times
  V\to \Q$
  is $\Q$-bilinear,
  non-degenerate
  and $\g\times G$-invariant. Write
  $\theta=\zeta_3-\zeta_3^2={\sqrt{-3}}$
  and define
  $B:V\times V \to \Q(\zeta_3)$ by
  $$B(u,v)=\theta^{-1}\left(\zeta_3 b(u,v)-
    b(u,\zeta_3 v)\right);$$
  then $B$ is a
  non-degenerate
  $G$-invariant Hermitian $\Q(\zeta_3)$-sesquilinear form
  and $b=\tr_{\Q(\zeta_3)/\Q}B$.
  
  Let us say that the signature of $B$ is $(p,r)$;
  then the signature of $b$ is $(2p,2r)$.
  We are dealing with elliptic surfaces of
  geometric genus $h$ and irregularity $q$, so that $p=h$
  and $r=4(h+1-q)$.
  So $p\ge 0$ and $r\ge 4$. In particular,
  $\dim_{\Q(\zeta_3)}V\ge 4$.

  Define a \emph{root} to be any $G$-conjugate of a basic
  cycle. Then, for every root $\delta$, $B(\delta,\delta)=-1$ and
  the corresponding complex reflexion $\s=\s_\delta$ is given by
  the formula
  $$\s_\delta(x)=x+(1-\zeta_3) B(x,\delta)\delta.$$
  This shows that $\s_\delta$ acts as a unitary reflexion of
  order $3$.
  The unitary group $U_B$ is an algebraic group
  defined over $\Q$ and has a determinant homomorphism
  $$\det :U_B\to R^1_{\Q(\zeta_3)/\Q}\GG_m$$
  that is also defined over $\Q$;
  here $R^1_{\Q(\zeta_3)/\Q}\GG_m$ is the kernel of the norm
  homomorphism $R_{\Q(\zeta_3)/\Q}\GG_m\to\GG_m$.
  The kernel of $\det$ is $SU_B$. The centre of $U_B$ is
  isomorphic to $R^1_{\Q(\zeta_3)/\Q}\GG_m$. Since each
  $\s_\delta$ lies in $U_B$ it follows that $G\subseteq U_B$
  and that therefore $M\subseteq U_B$.
  \begin{lemma}\label{perp}
    $V$ is an absolutely irreducible $\Q(\zeta_3)[G]$-module.
    Moreover, $V$ is spanned,
    as a $\Q(\zeta_3)$-vector space, by
    the $G$-orbit of any root.
    \begin{proof} Put
      $$V_\R=V\otimes_\Q\R=V\otimes_{\Q(\zeta_3)}\C.$$
      Note first that $(V_\R)^G
      =(V^G)\otimes_{\Q(\zeta_3)}\C=0$, by Lemma \ref{vanish},
      and that $B$ extends
      to a Hermitian form on $V_\R$.

      Suppose that $V'$ is a non-zero
      sub-$\C[G]$-module of $V_\R$.
      If there is a root $\delta$ such that
      $V'\subseteq Fix(\s_\delta)$, then
      $$V'=g(V')\subseteq Fix(g\s_\delta g^{-1})$$
      for all $g\in G$, so that,
      since the $G$-conjugates
      of $\s_\delta$ generate $G$,
      $V'\subseteq (V_\R)^G=0$.
      So every
      $\s_\delta$ acts non-trivially on $V'$.

      Fix a root $\delta$ and choose $x\in V'$ such that
      $\s_\delta(x)\ne x$.
      Then the formula defining $\s_\delta$ shows
      that $\delta\in V'$. Let $P$
      be the $\C$-subspace of $V_\R$ spanned by
      the cycles $\delta$; then $P$ is a
      non-zero $\C[G]$-module
      and $P$ is contained in every non-zero
      sub-$\C[G]$-module $V'$ of $V_\R$.
      In particular, $P$ is irreducible.

      Now suppose that $P^\perp\ne 0$. Then
      $P\subseteq P^\perp$. However, $\delta\in P$
      and $B(\delta,\delta)\ne 0$, so $P^\perp =0$.
      Therefore $P=V_\R$ and the lemma follows.
    \end{proof}
  \end{lemma}
  Now let $M$ denote the algebraic subgroup,
  over $\Q$, of the unitary group $U_B$
  that is generated by $G$. Since $G$ is generated by a
  single conjugacy class of unitary reflexions of order $3$,
  the same is true of $M$. Let
  $\tmu_n\subset R^1_{\Q(\zeta_3)/\Q}\GG_m$ be the $n$-torsion
  subgroup; this is the quadratic twist of $\mu_n$ determined by
  $\Q(\zeta_3)/\Q$.
  So $\det^{-1}(\tmu_3)$ is an extension of $\tmu_3$ by
  $SU_B$; we denote this group by $SU_B.\tmu_3$.
  \begin{theorem}\label{9.8}
    Either $M=SU_B.\tmu_3$
    or $G$ is finite.
    \begin{proof}  Observe first that
      $M\subseteq SU_B.\tmu_3$ since
      $G$ is generated by elements of order $3$.

      We set up some notation, as follows.
      \begin{enumerate}
      \item
        $\tM=R^1_{\Q(\zeta_3)/\Q}\GG_m.M\subset U_B$.
      \item
        $\Sigma\subset V_\R$
        is the real
        quadric hypersurface in $V_\R$
        defined by $B(x,x) =-1$.
      \item Take any root $\delta_0$
        and let $\tR$ denote its $\tM(\R)$-orbit
        and $R$ its $M(\R)$-orbit; these orbits
        are independent of the choice of $\delta_0$
        and are real semi-algebraic subsets of $\Sigma$.
      \item Given $\delta,\eta\in \tR$ that are not $\C$-proportional,
        let $L_{\delta,\eta}$ denote the complex $2$-plane
        that they span and let $H_{\delta,\eta}$ denote the real
        algebraic subgroup
        of $M_{\R}$ generated by $\s_\delta,\s_\eta$.
      \end{enumerate}
      Note that if $\delta\in\tR$ then $\delta=\lambda\eta$
      for some $\eta\in R$ and $\lambda\in\C$ with
      $\vert\lambda\vert =1$, so that we can define
      $\s_\delta$ by $\s_\delta=\s_\eta$. So $\s_\delta$
      is defined for all $\delta\in\tR$ and is a unitary
      reflexion of
      order $3$. Also $H_{\delta,\eta}$
      acts on $L_{\delta,\eta}$ as a group of unitary transformations
      and projects to a subgroup $P H_{\delta,\eta}$
      of the corresponding projective unitary group
      acting on $\P(L_{\delta,\eta}^\vee)=
      (L_{\delta,\eta}-\{0\})/\GG_m$.

      \begin{lemma}\label{linear} Assume that $\delta, \eta\in\tR$
        and are not $\C$-proportional.
        
        \part[i] If
        $\vert B(\delta,\eta)\vert >1$
        then $P H_{\delta,\eta}=PSU(1,1)\cong PSL_2(\R)$.
        \part[ii] If
        $\vert B(\delta,\eta)\vert <1$
        then, except for finitely many
        values of $\vert B(\delta,\eta)\vert$,
        $P H_{\delta,\eta}=PSU(0,2)\cong SO_3(\R)$.
        \part[iii]  If
        $\vert B(\delta,\eta)\vert =1$
        then $P H_{\delta,\eta}$ is a non-trivial
        split extension $\C\rtimes\Z/3$,
        regarded as a $2$-dimensional
        real algebraic group.
        \part[iv]
        There is a finite subset $\sC$ of the interval
        $[0,1)$ such that, if
        $\vert B(\delta,\eta)\vert\not\in \sC$, the group
        $P H_{\delta,\eta}$ is either $PSU(1,1)$ or
        $PSU(0,2)$ or  $\C\rtimes\Z/3$.
        \begin{proof} We start with the
          classification of the real algebraic subgroups
          $K$
          of $PSU_\b$, where $\b$ is a
          non-zero Hermitian
          form on a $2$-dimensional
          complex vector space of signature
          $(0,2),\ (1,1)$ or $(0,1)$, that
          are generated by elements of order $3$.
          In turn, this derives from the classification
          of the algebraic subgroups $H$ of
          the algebraic group $PSL_{2,\C}$ over $\C$.
          The fact that $3$ is odd simplifies matters.

          The algebraic subgroups $H$ of $PSL_{2,\C}$ are:
          $PSL_{2,\C}$ itself; subgroups of
          the normalizer of a maximal torus;
          subgroups of a Borel subgroup $B$;
          the
          polyhedral subgroups
          ${\mathfrak A}_4,{\mathfrak S}_4$ and
          ${\mathfrak A}_5$. So if $H$ is the Zariski
          closure of a group
          generated by elements of order $3$ then $H$ is one
          of: $PSL_{2,\C}$; $\Z/3$; a non-trivial
          split
          extension of $\Z/3$ by
          the unipotent radical $U\cong\C$ of $B$;
          ${\mathfrak A}_4$ or ${\mathfrak A}_5$.

          \begin{enumerate}
          \item
            It follows that any real algebraic subgroup $K$ of
            $PSU(1,1)\cong PSL_2(\R)$ that is the Zariski
            closure of a group generated by elements
            of order
            $3$ is either $\Z/3$ or $PSL_2(\R)$.
            The fact that
            $\delta$ and $\eta$ are not proportional
            shows that $P H_{\delta,\eta}\ne\Z/3$.
          \item
            The analogous
            subgroups of $PSU(0,2)\cong SO_3(\R)$ are $SO_3(\R)$,
            $\Z/3$, ${\mathfrak A}_4$ and ${\mathfrak A}_5$.
            The cases of
            ${\mathfrak A}_4$ and ${\mathfrak A}_5$ are excluded
            by excluding finitely many
            values of $\vert B(\delta,\eta)\vert$
            because we need to consider
            the situation
            where $\delta,\eta$ are vertices of a
            spherical triangle with angles $\pi/3,\pi/3$
            at $\delta,\eta$. The exceptional
            cases are those where the third angle
            is a rational multiple $a\pi/b$ of $\pi$,
            $1/3<a/b<1$ and
            $b\le 5$. It is clear that
            this can only happen for finitely many values
            of $\vert B(\delta,\eta)\vert$, all of which lie in
            $[0,1)$. As before, $\Z/3$ is excluded
            by the fact that $\delta, \eta$ are not proportional.
          \item In the degenerate case
            where $\vert B(\delta,\eta)\vert =1$, we see similarly
            that $PH_{\delta,\eta}$ is a non-trivial
            split extension of $U$ by $\Z/3$.
          \end{enumerate}
          Finally, \DHrefpart{iv} follows from
          \DHrefpart{i}--\DHrefpart{iii}.
        \end{proof}
      \end{lemma}

      If $\delta, \eta\in\tR$ are $\C$-proportional then
      $\s_\delta=\s_\eta^\pm$, so they generate
      a copy of $\Z/3$.
      \begin{lemma} Suppose that
        $\vert B(\delta,\eta)\vert\not\in \sC$
        and that $\delta,\eta$ are not $\C$-proportional.
        Define $O_{\delta,\eta}$ to be the
        $H_{\delta,\eta}$ orbit of $\delta$
        in the complex $2$-plane $L_{\delta,\eta}$. Then
        $L_{\delta,\eta}^0=\C^*.O_{\delta,\eta}$ is a dense
        semi-algebraic
        subset (in the real sense)
        of $L_{\delta,\eta}$
        and $L_{\delta,\eta}^0\cap \Sigma\subseteq \tR$.
        \begin{proof} This follows at once from Lemma \ref{linear}.
        \end{proof}
      \end{lemma}
      Suppose that $W\subseteq V_\R=V\otimes_{\Q(\zeta_3)}\C$ is a
      sub-$\C$-vector space containing
      a dense semi-algebraic subset
      $W^0$ such that
      \begin{enumerate}
      \item $W^0\cap\Sigma$
        is non-empty,
      \item $W^0\cap \Sigma$ is dense in $W\cap \Sigma$ and
      \item $W^0\cap\Sigma$ is contained in $\tR$.
      \end{enumerate}
      
      In particular, $W$ is spanned by $W\cap\tR$.

      For example, the one-dimensional
      such subspaces are exactly the
      lines $\C\delta$, where $\delta\in R$,
      and, provided that $\delta,\eta$ are not proportional
      and $\vert B(\delta,\eta)\vert\not\in \sC$ ,
      the complex $2$-plane $L_{\delta,\eta}$ is another.

    Assume until the end of the proof of
    Theorem \ref{9.8}
    that $W$ is maximal with respect to these properties.
    So certainly $W\ne 0$.

    \begin{proposition}\label{W1}
      Assume that $W\ne V_\R$. Then
      $G$ is finite.
      \begin{proof}
        Note first that $\tR\setminus W$ is not empty, since
        $\tR$ generates $V_\R$.

        We proceed to establish
        Lemmas \ref{2.10} and \ref{2.11}.

        \begin{lemma}\label{2.10}
          Suppose that $\delta\in\tR\setminus W$.
          Then the
          function $\b_\delta:W\to\R$
          defined by $\b_\delta(x)=\vert B(x,\delta)\vert$
          only takes values in $\sC$ when restricted
          to $W\cap\tR$.
          \begin{proof}
            Suppose that $\eta\in W\cap\tR$ (so that,
            in particular, $\eta$ is not $\C$-proportional to
            $\delta$) and that
            $\vert B(\delta,\eta)\vert\not\in \sC$.
            Then $L_{\delta,\eta}^0\cap\Sigma\subseteq \tR$. Now
            let $\eta$ vary over the subset $\sS_1$ of
            $W\cap\tR$ defined by the conditions that
            $\vert B(\delta,\eta)\vert\not\in \sC$
            and $\eta$ is not proportional to $\delta$.
            This is a semi-algebraic set.
            Put $\sS_2=\cup_{\eta\in U_1}L_{\delta,\eta}^0$,
            so that
            $$\sS_2\cap\Sigma\subseteq\tR.$$
            Moreover, $\sS_2$ is semi-algebraic;
            this can be seen
            by defining the semi-algebraic
            subset $\sS_3$ of $\sS_1\times M(\R)$ by
            $$\sS_3=\{(\eta,m) \vert \eta\in \sS_1\ \textrm{and}\ m\in H_{\delta,\eta}\}$$
            and observing that $\sS_2$
            is the image of $\C^*\times \sS_3$ under the map
            $$\C^*\times \sS_3 \to V\ :\ (z,\eta,m)
            \mapsto z.m(\eta).$$
            Since $L^0_{\delta,\eta}$ is dense in
            $L_{\delta,\eta}$, it follows that
            $\sS_2$ is dense in $\cup_{\eta\in U_1}L_{\delta,\eta}$
            and so is dense in $W\oplus\C\delta$.
            Therefore
            we can enlarge $W$ to $W\oplus\C\delta$
            and then take $(W\oplus\C\delta)^0=\sS_2$ to
            contradict the maximality of $W$.
          \end{proof}
        \end{lemma}
        \begin{lemma}\label{2.11}
          \part[i] Suppose that
          $\delta\in\tR\setminus W$.
          Then \emph{either}
          $\delta\in W^\perp$
          \emph{or}
          $W\cap W^\perp$ is of $\C$-codimension
          one in $W$. In both cases $\delta$ is orthogonal
          to $W\cap W^\perp$.

          \part[ii]
          $W\cap W^\perp=0$.
          \begin{proof}\DHrefpart{i}
            $$W\cap\tR\supseteq W\cap W^0\cap\Sigma =W^0\cap\Sigma,$$
            so that, by Lemma \ref{2.10},
            $\b_\delta$ takes values only in the finite set
            $\sC$ on the non-empty
            dense semi-algebraic subset
            $W^0\cap\Sigma$ of $W\cap\Sigma$.
            Therefore $\b_\delta$ takes values only in $\sC$
            on the real quadric hypersurface
            $W\cap\Sigma$.

            Fix $c\in \sC$ and consider the condition that
            $\b_\delta(w)=c$
            for $w\in W$.
            There are complex co-ordinates $(z_k)_{k\in K}$ on $W$
            and disjoint subsets $I,J$ of $K$ such that
            $J$ is not empty,
            $W\cap\Sigma$ is given by the equation $D=0$,
            where
            $$D=\sum_{i\in I}\vert z_i\vert^2 -
            \sum_{j\in J}\vert z_j\vert^2+1,$$
            and the condition that $\b_\delta(w)=c$ is given by the equation
            $E=0$, where
            $$E=\vert B(z,\delta)\vert^2-c^2.$$
          Then $dD$ is proportional to $dE$ on $W\cap\Sigma$.
          Say $dE=\phi dD$ on $W\cap\Sigma$
          for some function $\phi$ on $W\cap\Sigma$.
          Write  $B(w,\delta)=\sum_{i\in I} a_iz_i-\sum_{j\in J} a_jz_j
          +\sum_{k\not\in I\cup J} b_k z_k.$
          Then
          $$\bar{a}_k B(w,\delta)=\phi z_k$$
          on $W\cap \Sigma$
          for all $k\in I\cup J$.
          This gives two possibilities: \emph{either}
          {
            \renewcommand{\labelenumi}{(\alph{enumi})}
            \begin{enumerate}
            \item
              the function $w\mapsto B(w,\delta)$ is identically zero
              on $W\cap\Sigma$ (in which case it vanishes
              on $W$ and $\delta\in W^\perp$) \emph{or}
            \item
              $I$ is empty,
              $J$ has just one element,
              $W\cap\Sigma$ is the cone $\widehat{\Sigma}$
              over
              a circle $\vert z\vert^2=1$
              and the restriction of $\b_\delta$ to $W\cap\Sigma$
              is proportional to $\vert z\vert$.
            \end{enumerate}
            }
            \noindent In the first case $\delta\in W^\perp$.
            in the second case $\b_\delta$
            vanishes on the vertex
            of $\widehat{\Sigma}$ and this vertex has codimension
            $1$ in $W$. Moreover, this vertex is $W\cap W^\perp$,
            and part $(1)$ of the lemma is
            proved.
            \DHrefpart{ii}
            Note that $\delta$ is orthogonal
            to $W\cap W^\perp$, from $(1)$.
            Then $W\cap W^\perp$ is orthogonal
            to every root in $W$ and to every root not in $W$,
            so lies in $V_\R^G$. But this space vanishes,
            by Lemma \ref{vanish}, and the lemma is proved.
          \end{proof}
        \end{lemma}

        Now we finish the proof of Proposition
        \ref{W1}.
        
        First, suppose that $\dim_\C W\ge 2$. Then,
        by Lemma \ref{2.11}, $\tR\setminus W\subseteq W^\perp$,
        so that for every root in $\tR\setminus W$
        the reflexion $\s_\delta$ preserves $W$.
        Certainly $\s_\delta$ preserves $W$ for every
        root $\delta$ in $W$, so that $W$ is $G$-invariant.
        But $V_\R$ is irreducible and $W\ne V_\R$, so
        we have a contradiction and $\dim_\C W=1$.
        
        Next, since $W\cap\tR$ is not empty, $W\cap R$ is also non-empty
        and then $W$ is a line $\C\delta$ for some $\delta\in R$.
        
        Next, let $\tM^0$ denote the identity
        connected component of $\tM$
        and $\tR^0=\tM^0(\R)\delta_0$.
        Then, for all $\delta\in\tR^0$
        and for all $m\in\tM^0(\R)$ close to
        the identity,
        $\vert B(\delta, m(\delta))\vert$
        is close to $1$ and so does not lie in $\sC$.
        Therefore $\delta$ and $m(\delta)$
        are
        proportional, since otherwise $L_{\delta,\eta}$
        would be a $2$-dimensional space $W$. So $\tM^0$
        preserves each
        complex line $W$.
        Since the lines $\C\delta$ span $V_\R$,
        it follows that $\tM^0$ is a real algebraic torus
        and that $V_\R$ decomposes as a
        direct sum
        $V_\R=\oplus U_\chi$,
        where $\chi$ runs over the complex characters
        of $\tM^0$ and $U_\chi$ is the
        $\chi$-eigenspace of $V_\R$.

        Finally, suppose that $r$ is the number of distinct
        characters in this decomposition
        and that $r\ge 2$.
        Since $\tM^0$ is normal
        in $\tM$, this decomposition is
        a system of imprimitivity for
        the action of $\tM$ on $V_\R$,
        and so for the group $G$.
        However, if $G$ were imprimitive
        then it would, since the representation
        $V_\R$ is irreducible, possess
        a surjection onto the symmetric group
        ${\frak S}_r$ (\cite{ST}, penultimate
        paragraph of
        p. 276; at this point of their paper
        the assumptions that the group
        in question be finite and unitary are not relevant).
        Since $G$ is generated
        by elements of odd order
        this is impossible, so that
        $r=1$ and $\tM^0$ is $1$-dimensional. Since $\tM^0$ contains $S^1$
        it is then equal to $S^1$, which is
        the centre of $U_B$. But
        $M\subseteq SU_B.\tmu_3$,
        so $M$ is finite
        and Proposition \ref{W1}
        is proved.
      \end{proof}
    \end{proposition}

    \begin{proposition} Assume that
      $W=V_\R$. Then $M=SU_B.\tmu_3$.
      \begin{proof} Since $W=V_\R$, it follows that
        $V_\R^0\cap\Sigma\subseteq\tR\subseteq\Sigma$,
        so that $\tR$ is dense in $\Sigma$.
        Therefore $M(\R)$ equals the
        real algebraic subgroup
        $N$ of $U_B(\R)$ that is generated by the
        unitary reflexions of the form
        $$x\mapsto x+(1-\zeta_3)B(x,\delta)\delta,$$
        as $\delta$ runs over the elements of
        $\Sigma$. It is clear that $N$
        is normal in $U_B(\R)$ and is contained in
        $(SU_B.\tmu_3)(\R)$. Since $SU_B(\R)$ is a simple
        real algebraic group, it follows that
        $M(\R) =(SU_B.\tmu_3)(\R)$ and
        the proposition is proved.
      \end{proof}
    \end{proposition}
    This completes the proof of Theorem \ref{9.8}.
  \end{proof}
\end{theorem}
\begin{proposition} If $G$ is finite then
  $h=q=0$ and $G$ is the group
  $G=W(L_4)$ (the group
  numbered
  $32$
  in the list constructed by
  Shephard and Todd \cite{ST}.)
  \begin{proof} We know that $G$ is an irreducible finite
    complex reflexion group.
    So it
    is one of the groups in the Shephard--Todd list.
    The fact that it is generated by reflexions
    of order $3$ leaves only a few possibilities.
    
    Let $n$ denote the dimension of its defining
    representation, which is $V_\R$,
    so that
    $$n=\dim_{\Q(\zeta_3)}V_\R=
    h^{1,1}_{prim}(X)/2=5h+4(1-q)=h+4(h-q+1).$$
    Since $h\ge q\ge 0$ it follows that
    $n\ge h+4\ge 4$. Say that $G=ST_N$, the group
    (or family of groups) numbered $N$ in the list of \cite{ST}
    and consider the table on p. 412 of \cite{Coh76},
    which covers the cases $N\in [24, 34]$.
    Taking into account the facts that $G$ is generated
    by reflexions of order $3$ and that $n\ge 4$ shows that,
    if $N\in [24, 34]$, then
    $N=32$ and $n=4$, so that $h=0$
    and then $q=0$.
    For $N\le 23$ look at the table on p.
    $301$ of \cite{ST}; since $n\ge 4$ we have
    $N=1$ or $2$. However, each group in either of these families
    admits a surjection onto a non-trivial symmetric group,
    so is not generated by elements of order $3$.
    If $N\ge 35$ then $G$ is a Weyl group of type $E$,
    which groups are not generated by elements of order $3$.
  \end{proof}
\end{proposition}
\begin{remark}
  When $h= q =0$
  then
  $X$ is the blow-up of a del Pezzo
  surface $S$ of degree $1$ at the base point of
  $\vert -K_S\vert$. This shows that the
  case where $G=ST_{32}$ does occur, and echoes the existence
  of an embedding $W(L_4)\inj W(E_8)$.
\end{remark}
We continue with a fixed generic curve $C$ of genus $q\ge 2$
and a
generic line bundle $L$ on
it of degree $h+1-q\ge 1$ and consider the Eisenstein surfaces $X_6$
defined by
these data. For a Jacobian elliptic surface $Y$
we denote by ${\HS}(Y)$ the $\Z$-Hodge structure on
$H^2_{prim}(Y,\Z)$. The negative of the cup product defines
a polarization on this; we call it the standard polarization.

\begin{theorem}\label{generic} Suppose that $f:X_6\to C$ is a
  very general Eisenstein surface and that the assumption
  (\ref{assume}) holds. Then the natural
  homomorphism $\phi$ from $\Aut(X_6)$
  to the automorphism group $\Aut(\HS(X_6))$
  of the polarized
  $\Z$-Hodge structure $\HS(X_6)$
  is an isomorphism.
  \begin{proof} Suppose that $1\ne \g\in\ker\phi$. Then
    $\g$ acts trivially on $\vert K_{X_6}\vert$.
    Since $K_{X_6}\sim f^*(K_C+M)$, the assumption
    (\ref{assume}) ensures that $\vert K_{X_6}\vert$
    pulls back from a very ample class on $C$,
    so that $\g$ acts trivially on $C$ and the geometric
    quotient $Y=[X_6/\langle\g\rangle]$ admits a morphism
    to $C$ whose generic fibre is the quotient of
    the generic fibre of $f$ by $\langle\g\rangle$.
    So $Y$ is birationally ruled over $C$, so that
    $p_g(Y)=0$. On the other hand
    $H^0(X_6,\Omega_{X_6}^2)=H^0(X_6,\Omega_{X_6}^2)^\g$,
    so vanishes,
    which is absurd. So $\phi$ is injective.

    To prove the surjectivity of $\phi$ we argue as follows.
    By \cite{Del72}, Proposition 7.5,
    the monodromy group $G$ has a subgroup $G_0$ of finite index
    that embeds
    as a subgroup of the Mumford--Tate group $MT(X_6)$
    of the polarized $\Z$-Hodge structure of a very
    general surface $X_6$ in the family of Eisenstein surfaces
    defined by $(C,L)$.

    \begin{enumerate}
    \item $MT(X_6)$
      is a $\Q$-subgroup of $O_b$ but, because of
      the existence of the automorphism $\g$,
      $MT(X_6)$ is in fact a $\Q$-subgroup of
      $U_B$.
    \item The identity connected component $M^0$
      of $M$ is contained in $MT(X_6)$
      since $M^0$ is contained
      in the $\Q$-algebraic subgroup generated by
      $G_0$.
    \item So $MT(X_6)$ contains $SU_B$.
      Since $MT(X_6)$ contains the scalars,
      it equals $U_B$ and its centralizer is $\Q(\zeta_3)$.
    \end{enumerate}

    If $\theta$ is an automorphism of the polarized
    $\Z$-Hodge structure
    $\HS(X_6)=V_\Z$ then $\theta$
    normalizes the image of the Deligne torus
    $\mathbb{S}=R_{\C/\R}\GG_m$ in $GL(V_\R)$ that
    defines $\HS(X_6)_\R$. Examination of the Galois module
    $\X^*({\mathbb{S}})$ shows that the $\R$-automorphism
    group of $\mathbb{S}$ consists of the
    four maps
    $$z\mapsto\{z,z^{-1},{\overline{z}},{\overline{z}}^{-1}\}.$$
    Since the weight of $\HS(X_6)$ is non-zero,
    $\theta$ acts trivially on the copy of $\R^*$ in
    $\mathbb{S}$ and so its action on $\mathbb{S}$
    is either $z\mapsto z$ or $z\mapsto {\overline{z}}$. Since $H^{2,0}\ne 0$
    it follows that $\theta$ acts trivially on $\mathbb{S}$.

    Let $Z$ denote the connected
    component of the fixed locus of $\theta$ in $GL(V)$;
    this is a $\Q$-algebraic subgroup of $GL(V)$ and
    $Z_\R$ contains $\mathbb{S}$. So $Z$ contains
    $MT(X_6)$, so that $\theta$ centralizes $MT(X_6)$ and lies in
    $\Q(\zeta_3)$. But
    also
    $\theta$ and $\theta^{-1}$ are integral over $\Z$,
    since they act on the finite $\Z$-module $\HS(X_6)$, so
    that they both lie in $\Z[\zeta_3]$
    and then $\theta$ is a power of $\zeta_6$.
  \end{proof}
\end{theorem}
\begin{corollary}\label{2.17}
  (Generic local Torelli) If $X_6$ is a generic Eisenstein surface
  and $\vert K_{X_6}\vert$ has no base points then
  the coarse period map
  $$per:[\sje]\to [\sD/\G]$$
  is injective in an analytic \nbd of the point corresponding to $X_6$.
  \begin{proof} It's enough to allow $X_6$ to be very general.

    Since $\vert K_{X_6}\vert$ has no base points,
    the infinitesimal Torelli theorem holds for $X_6$.
    The result follows at once from the theorem.
  \end{proof}
  \end{corollary}
\begin{remark} Really the important point here is the weaker
  result that the automorphism group of the polarized Hodge
  structure $\HS$ is induced by the automorphism group of the
  variety
  modulo automorphisms of $\HS$ that act trivially on the period
  domain.
  So very general Eisenstein surfaces happen to behave like hyperelliptic
  curves rather than non-hyperelliptic ones.
\end{remark}
\end{section}
\begin{section}{$\Q$-Hodge structures and generic Torelli}
  Suppose that $n$ is an integer. If $X,Y$ are Jacobian elliptic
  surfaces then
  an $n$-isogeny from $\HS(X)$ to $\HS(Y)$ is a $\Z[1/n]$-isomorphism
  $\psi:\HS(X)\otimes\Z[1/n]\to \HS(Y)\otimes \Z[1/n]$ of
  $\Z[1/n]$-Hodge structures
  such that
  $n\psi(\HS(X))\subseteq \HS(Y)\subseteq n^{-1}\psi(\HS(X))$
  and $\psi$ preserves the standard polarizations
  up to a positive rational scalar multiple.
  This notion defines a
  correspondence
  $\G_n$ on $\sD/\G$ and on $[\sD/\G]$ and then we can consider
  the composite
  $$per_n=\G_n\circ per:[\sje_{h,q}]\to[\sD/\G]\to [\sD/\G].$$
  This is multi-valued; we say that it is \emph{very generally
    injective}
  if, for two points $X,Y$ in $[\sje_{h,q}]$ of which one is very
  general, we have
  $X=Y$ if the finite set $per_n(X)$ has non-empty intersection
  with $per_n(Y)$.
  \begin{lemma}\label{more generic}
    If $\HS(Y)$ is isogenous to $\HS(X_6)$ for a
    very general Eisenstein surface $X_6$ then the automorphism
    group of $\HS(Y)$, with its polarization, is either
    $\langle\zeta_6\rangle$ or $\{\pm 1\}$
    \begin{proof} $\HS(Y)_\Q$ is isomorphic to $\HS(X)_\Q$.
      The proof of Theorem \ref{generic} then shows, first, that every
      automorphism $\theta$ of $\HS(Y)$ lies in $\Q(\zeta_3)$
      and, second, that $\theta$ is a power of $\zeta_6$.
      Since $-1$ is always an automorphism of $\HS(Y)$
      the lemma follows.
    \end{proof}
  \end{lemma}
  \begin{theorem} Under the assumptions (\ref{assume})
          the multi-valued map $per_n$ is very generally injective.
          \begin{proof} Assume otherwise. As recalled in Section
            \ref{intro} the IVHS determines the base curve
            $C$ and the ramification divisor $Z$ of the classifying
            morphism $\Phi:C\to\bsell$, and then for generic
            $\Phi$ the pair $(C,Z)$ determines $\Phi$
            modulo the action of the automorphism
            group $\GG_m$ of $\bsell$. Since the closure of such an
            orbit is a Gauss--Eisenstein pencil, a generic point $X$
            of $\sje$ determines a Gauss--Eisenstein pencil, the
            closure            of the $\GG_m$-orbit through the
            point $X$, and any
            failure of $per_n$ to be very generally
            injective
            can be detected on a very general Gauss--Eisenstein pencil.

            Suppose that $X_6$ is a very general Eisenstein surface
            and that
            $\sX\to B \cong \P^1$ is the Gauss--Eisenstein pencil
            generated by $X_6$ and some Gauss surface $X_4$. Let
            $P\in B$ denote the point corresponding to $X_6$,
            $B'$ the normalization of
            the image of $B$ in the coarse moduli space $[\sje^{RDP}_{C,L}]$
            and $P'\in B'$ the image of $P$. In an analytic neighbourhood
            of $P$ the map $B\to B'$ is just the geometric quotient
            by the group $\langle\zeta_6\rangle$. (This action
            is not effective; the kernel is
            $\{\pm 1\}$).
            
            Consider the multi-valued map
            $per_n\vert_{B'}:B'\to [\sD/\Gamma]$; by assumption, this is
            not very generally injective in the sense defined above.
            However, Corollary \ref{2.17} and Lemma \ref{more generic}
            combine to show that
            $per_n\vert_{B'}$ is an isomorphism onto each of its
            images $Y$
            \emph{either}
            \begin{enumerate}
              \item in an analytic \nbd
            of $P'$ in $B'$ (when
            $\Aut(\HS(Y))=\langle\zeta_6\rangle$) \emph{or}
            \item in an analytic \nbd
              of $P$ in $B$ (when $\Aut(\HS(Y))=\{\pm 1\}$).
            \end{enumerate}
            Since $B'$ is complete there is,
            therefore, another point $Q'\in B'$
            such that $Q'\ne P'$ and $per_n(P')$ meets $per_n(Q')$.
            Let $Q$ be a pre-image of $Q'$ in $B$.
            Varying the pencil $B$ by
            keeping the Eisenstein
            surface $X_6$ fixed and varying the Gauss
            surface $X_4$
            then gives
            a locus $\sL$ in $\sje^{RDP}_{C,L}$,
            swept out by the points $Q$,
            that consists of
            elliptic surfaces whose Hodge structure
            is $n$-isogenous to that of $X_6$.
            Since a general member of $\sje$ determines
            the Gauss--Eisenstein pencil on which it lies,
            the point $Q$ varies when $X_4$ varies,
            so that $\dim\sL$ is strictly positive.
            However,
            this contradicts the infinitesimal Torelli theorem.
          \end{proof}
        \end{theorem}
        \begin{corollary} If $8h>10(q-1)$, $h\ge q+3$ and $q\ge 2$
          then a very general surface in $\sje_{h,q}$ is
          determined by its polarized $\Q$-Hodge structure.
          \noproof
        \end{corollary}
        
              \end{section}

\end{document}